 \documentclass[final,3p,times]{elsarticle}
\usepackage{amsmath,hyperref}
\usepackage{lineno}

\usepackage{epstopdf}
\journal{Elsevier}

\usepackage{lineno,hyperref}
\usepackage{graphicx}% Include figure files
\usepackage{dcolumn}% Align table columns on decimal point
\usepackage{bm}% bold math
\usepackage{epsfig}
\usepackage{booktabs}
\usepackage{subfigure}
\usepackage{graphics}
\usepackage{amssymb}
\usepackage{amsmath}
\usepackage{array}
\usepackage{color}
\usepackage{booktabs}
\usepackage{multirow}
\usepackage{caption}
\usepackage{chngpage}
\usepackage{subfigure}
\usepackage{mathrsfs,gensymb,float}

\biboptions{numbers,sort&compress}
\modulolinenumbers[1]

\begin{document}

%\linenumbers   %显示行数
\captionsetup[figure]{labelfont={bf},name={Fig.},labelsep=period}        %改变图编号为Fig.

\begin{frontmatter}
	
\title{An efficient thermal lattice Boltzmann method for simulating three-dimensional liquid-vapor phase change}

\author[mymainaddress,myaddress2]{Jiangxu Huang}
\author[mymainaddress,myaddress2]{Lei Wang\corref{mycorrespondingauthor}}
\ead{wangleir1989@126.com}
\author[mymainaddress,myaddress2]{Kun He}
\author[mymainaddress,myaddress2]{Changsheng Huang}

\cortext[mycorrespondingauthor]{Corresponding author}

\address[mymainaddress]{School of Mathematics and Physics, China University of Geosciences, Wuhan 430074, China}
\address[myaddress2]{Center for Mathematical Sciences, China University of Geosciences, Wuhan 430074, China}

\begin{abstract}
In this paper, a multiple-relaxation-time lattice Boltzmann (LB) approach is developed for the simulation of three-dimensional (3D) liquid-vapor phase change based on the pseudopotential model. In contrast to some existing 3D thermal LB models for liquid-vapor phase change, the present approach has two advantages: for one thing, the current approach does not require calculating the gradient of volumetric heat capacity [i.e., $\nabla \left( {\rho {c_v}} \right)$], and for another, the current approach is constructed based on the seven discrete velocities in three dimensions (D3Q7), making the current thermal LB model more efficient and easy to implement. Also, based on the scheme proposed by Zhou and He [Phys Fluids 9:1591-1598, 1997], a pressure boundary condition for the D3Q19 lattice is proposed to model the multiphase flow in open systems. The current method is then validated by considering the temperature distribution in a 3D saturated liquid-vapor system, the $d^2$ law and the droplet evaporation on a heated surface. It is observed that the numerical results fit well with the analytical solutions, the results of the finite difference method and the experimental data. Our numerical results indicate that the present approach is reliable and efficient in dealing with the 3D liquid-vapor phase change.     
\end{abstract}


\end{frontmatter}

\section{Introduction}
Liquid-vapor phase change is ubiquitous in nature and is often encountered in many areas of applied science and engineering technologies, such as thin thermal coating \cite{WanJoule2018}, chip manufacturing \cite{DugasLangmuir2005,MarcinichenAE2012}, and food preservation \cite{AdibJFE2008}. The importance of liquid-vapor phase change in modern industry is now well established, however, due to the complex physical effects at the liquid-vapor interface, it is still challenging to simulate liquid-vapor flows with phase change \cite{TaylorIJHMT2009,DadhichRSER2010,KakacIJHMT2008}.  

As a mesoscopic numerical approach, in the past two decades, the lattice Boltzmann method (LBM) has become a numerically robust and efficient technique for simulating complex fluids \cite{AidunARFM2010,Kruger}, particularly multiphase flows \cite{HuangJWS2015,LiPECS2016,HuangPRL2021,LiIJMF2022} and porous media flows \cite{HeIJHMT2019}. As for liquid-vapor flows with phase change, researchers in the LB community have proposed several phase change models and the most common ones are the phase-field-based \cite{DongNHTPA2009,SunCF2013,SafariPRE2010} and the pseudopotential-based models \cite{ZhangPRE2003,MarkusPRE2011,GongIJHMT2012,LiPRE2017,ZhangPRE2021}. In the phase-field-based phase change model, owing to the source term added in the governing equations is related to the spatial gradient of the order parameter, one must add a small perturbation of the liquid-vapor interface or a small vapor bubble at the initial time to trigger the liquid-vapor phase change \cite{DongNHTPA2009}. In such a case, the phase-field based phase change model is unable to model the onset of nucleation boiling. 

Unlike the phase-field-based model, the flow and temperature fields in the pseudopotential-based model are coupled via the equation of state \cite{YuanPOF2006}. In such a case, the non-ideal pressure drives the liquid-vapor phase change, and then the phase interface can spontaneously arise and deform \cite{HuangJWS2015,LiPECS2016,Kruger}. In 2003, Zhang and Chen \cite{ZhangPRE2003} first proposed a phase change model using the van der Waals type of gas model to couple the energy transport equation and the pseudopotential model. Later, Hazi and Markus \cite{MarkusPRE2011} modified the energy equation according to the local balance law for entropy and proposed a thermal pseudopotential LB model for liquid-vapor phase change. However, the equation of state (EOS) employed in Hazi and Markus'model is an artificial EOS, not an EOS for real gases. Gong and Cheng \cite{GongIJHMT2012} then developed another thermal LB model to simulate liquid-vapor phase change heat transfer in this setting, where the Peng-Robinson EOS is adopted. After that, Li et al. \cite{LiPRE2017} pointed out that the assumption used for the thermal diffusion term in Gong and Cheng's model is unreasonable in multiphase flows, and they further developed an improved therm LB phase change model. Recently Zhang \cite{ZhangPRE2021} also proposed a thermal LB model for liquid vapor phase change, where the latent heat of vaporization is decoupled with the EOS, and the calculation of the Laplacian term that appeared in Li et al.'s model is avoided. Although some thermal LB models have been constructed to simulate liquid-vapor flows with phase change, there exist two primary defects in these approaches. To begin with, to recover the correct energy equation,  the source terms used in these models contain the gradient term of $\nabla \left( {\rho {c_v}} \right)$ (here,  $ \rho $ is the fluid density and $c_v$ is the specific heat at constant volume), which may bring some additional errors when the fluid properties across the liquid-vapor interface are significantly different \cite{HuATE2019,HuCF2019}, and meanwhile, this treatment also increases the computational load in the practical simulations, and this phenomenon is more distinct for 3D cases. Secondly, all the mentioned thermal LB models are constructed for the two-dimensional (2D) cases. They can not be directly used for the 3D liquid-vapor phase change problems due to the significant differences in the development and implementation of 3D LB models. These above issues may also be why the finite-difference method (FDM) is still the most widely used approach in LB community for modeling 3D liquid-vapor phase change today \cite{FeiPOF2019,QinJFM2019,LiATE2020,LuoPTRSA2021,CaiICHMT2022}. More recently, Fogliatto et al. \cite{FogliattoIJTS2022} developed a multiple-relaxation-time (MRT) LB model for boiling heat transfer with the assumption of the uniform distribution of the specific heat at constant volume. Although the temperature equation can be recovered correctly with the proposed model, it is constructed based on the  D3Q15 lattice, which is not efficient for 3D simulations. Moreover, just like previous 2D approaches, this model needs to calculate the density gradient, and we note that this issue also exists in another 3D thermal LB model proposed by Li et al \cite{LiPRE2022}.  

To address the above problems, in the current work,  we target to construct an efficient thermal LB model for simulating 3D liquid-vapor phase change, and it can be viewed as an extension of our recent 2D LB approach for non-ideal fluids \cite{WangPRE2022}. As shown later, the aforementioned defects do not exist for the present model such that it holds the major advantages of the LB method. Furthermore, we also develop a constant-pressure boundary condition on D3Q19 lattice to model the liquid-vapor phase change in open system. The remainder of the present paper is organized as follows. In Section 2, the pseudopotential multiphase model together with the proposed thermal LB model are presented. In Section 3, the frequently used boundary condition for liquid-vapor phase change is introduced. Validations and discussions are given in Section 4, and a brief conclusion is drawn in Section 5.

\section{Lattice Boltzmann model}
\begin{figure}[H]
	\centering
	\subfigure[]{ \label{fig1a}
		\includegraphics[width=0.5\textwidth]{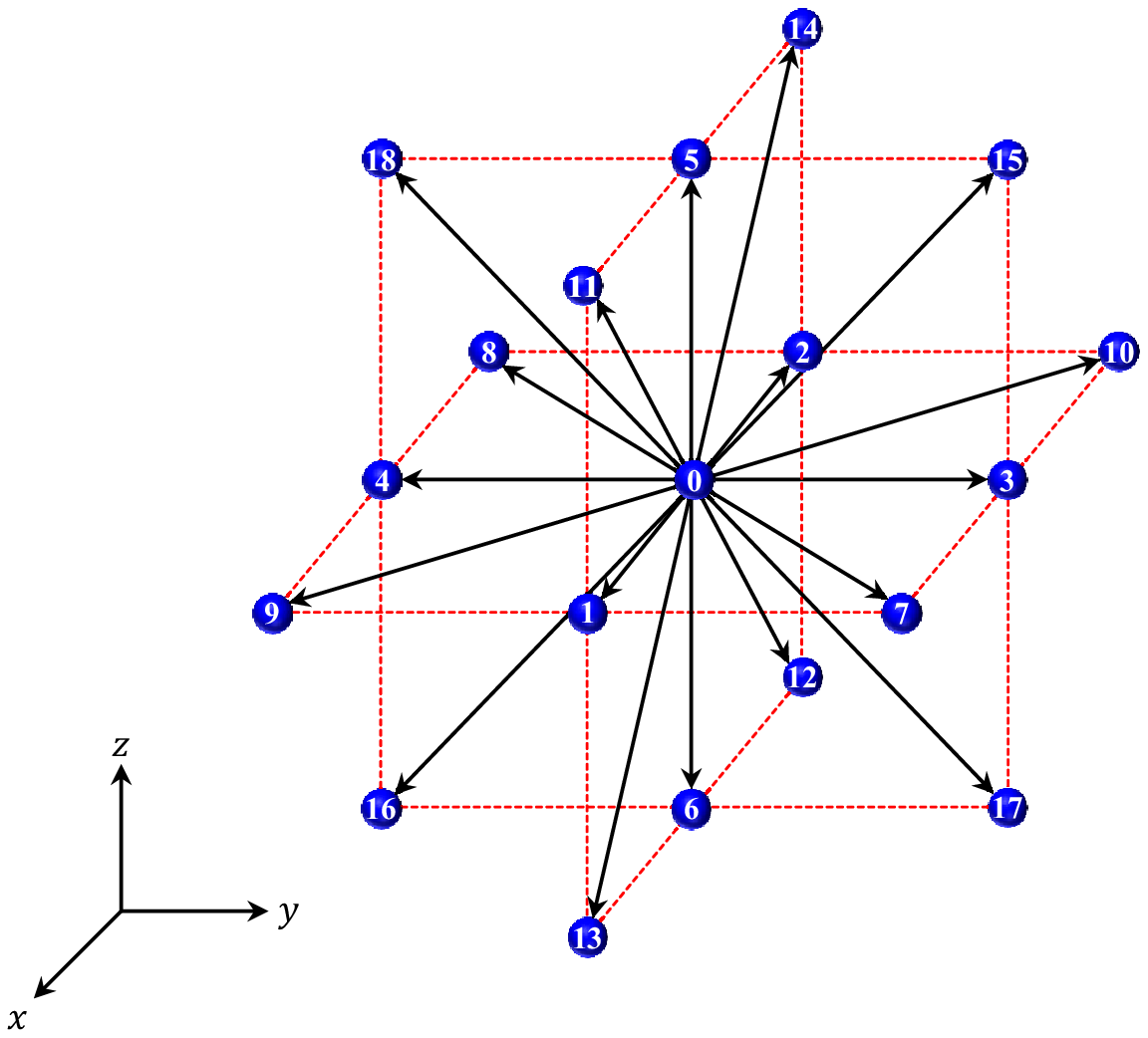}}
	\subfigure[]{ \label{fig1b}
		\includegraphics[width=0.375\textwidth]{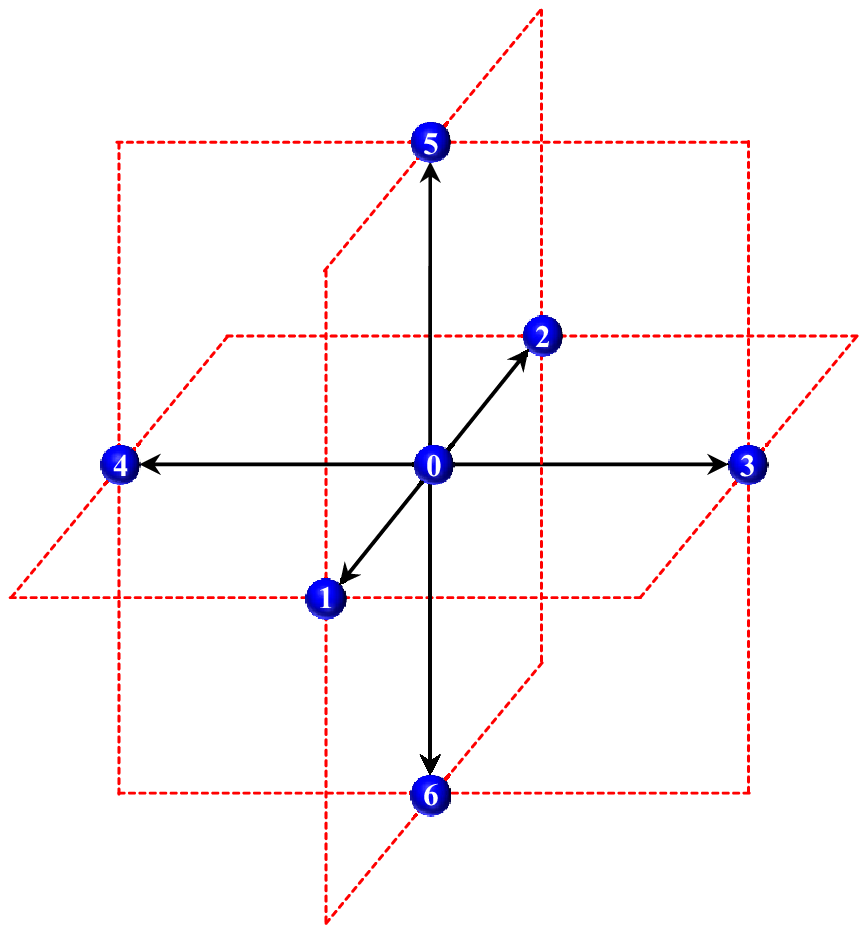}}\\		
	\caption{D3Q19 (a) and D3Q7 (b) velocity sets.    }
	\label{fig1}
\end{figure}

The present thermal model for liquid-vapor phase change is constructed based on the pseudopotential multiphase LB model \cite{ShanPRE1993}. For this approach, the liquid-vapor phase change is driven by the equation of state $p_{EOS}$ without using any artificial phase-change terms, and it adopts two distribution functions: one is the density distribution function for the velocity field based on the D3Q19 lattice \cite{LiCF2019}, and the other is the thermal distribution function for the temperature field based on the D3Q7 lattice.  The discrete velocity sets of the above two lattice models are schematized in Fig. \ref{fig1}. 

As is well known that the models in LB community can be classified into three kinds on the basis of the collision operator: the Bhatnagar-Gross-Krook (BGK) model \cite{QianEPL1992}, the two-relaxation-time (TRT) model \cite{GinzburgCCP2008} and the MRT model \cite{LuoMRT2002}. In this paper, we adopt the MRT model for its better numerical accuracy and stability. Moreover, the transformation matrix used here is a non-orthogonal one, which is much simpler than the orthogonal matrix adopted in most previous works \cite{HeIJHMT2019}.

\subsection{Lattice Boltzmann equation for density function}
The evolution equation of the non-orthogonal MRT model for velocity filed can be written as \cite{LiCF2019} 
\begin{equation}
{f_i}\left( {{\bf{x}} + {{\bf{c}}_i}\Delta t,t + \Delta t} \right) - {f_i}\left( {{\bf{x}},t} \right) = {\left( {{{{\bf{\hat M}}}^{ - 1}}{\bf{\hat \Lambda \hat M}}} \right)_{ij}}\left[ {{f_i}\left( {{\bf{x}},t} \right) - f_i^{\left( {eq} \right)}\left( {{\bf{x}},t} \right)} \right] + {{\bf{F}}_i}\left( {{\bf{x}},t} \right),
\end{equation}
where ${f_i}\left( {{\bf{x}},t} \right)$ is the discrete velocity distribution function at position $\bf{x}$ and $t$, ${{{\bf{c}}_i}}$ is the discrete velocity in $i$ direction (see Fig. 1), ${f_i^{\left( {eq} \right)}\left( {{\bf{x}},t} \right)}$ is the  equilibrium distribution function defined as \cite{ChenARFM1998}
\begin{equation}
f_i^{\left( {eq} \right)}\left( {{\bf{x}},t} \right) = {{\hat w}_i}\rho \left[ {1 + \frac{{{{\bf{u}}^{eq}} \cdot {{\bf{c}}_i}}}{{\hat c_s^2}} + \frac{{{{\left( {{{\bf{u}}^{eq}} \cdot {{\bf{c}}_i}} \right)}^2}}}{{2\hat c_s^4}} - \frac{{\left| {{{\bf{u}}^{eq}}} \right|}}{{2\hat c_s^2}}} \right],
\end{equation}
where ${{\hat c}_s}$ is the lattice sound speed and it is given by ${{\hat c}_s} = {{\Delta x} \mathord{\left/
{\vphantom {{\Delta x} {\left( {\Delta t\sqrt 3 } \right)}}} \right.
\kern-\nulldelimiterspace} {\left( {\Delta t\sqrt 3 } \right)}}$ with ${\Delta x}$ and ${\Delta t}$ being the lattice spacing and time step, respectively, and the values of them are both set to 1 in this work.  ${{\hat w}_i}$ is the weight associated with the velocity ${{{\bf{c}}_i}}$, which is defined as ${{\hat w}_0} = {1 \mathord{\left/
{\vphantom {1 3}} \right.
\kern-\nulldelimiterspace} 3},{{\hat w}_{1 - 6}} = {1 \mathord{\left/
{\vphantom {1 {18}}} \right.
\kern-\nulldelimiterspace} {18}},{{\hat w}_{7 - 18}} = {1 \mathord{\left/
{\vphantom {1 {36}}} \right.
\kern-\nulldelimiterspace} {36}}$ for the D3Q19 lattice. ${{\bf{\hat M}}}$ is the transformation matrix shown in Appendix. $\rho $ and ${{{\bf{u}}^{eq}}}$ are the fluid density and the equilibrium velocity defined as 
\begin{equation}
	\rho  = \sum\limits_{i = 0}^{18} {{f_i}} ,\;\;\;\;\;\rho {{\bf{u}}^{eq}} = \sum\limits_{i = 0}^{18} {{{\bf{c}}_i}{f_i}} .
\end{equation}
${{\bf{F}}_i}\left( {{\bf{x}},t} \right)$ is the exact-difference-method forcing term given by Kupershtokh et al. \cite{KupershtokhCMA2009}
\begin{equation}
{{\bf{F}}_i}\left( {{\bf{x}},t} \right) = f_i^{\left( {eq} \right)}\left( {{{\bf{u}}^{eq}} + \Delta {\bf{u}}} \right) - f_i^{\left( {eq} \right)}\left( {{{\bf{u}}^{eq}}} \right),
\end{equation} 
where $\Delta {\bf{u}}\left( {{\bf{x}},t} \right) = {{{\bf{F}}\left( {{\bf{x}},t} \right)\Delta t} \mathord{\left/
{\vphantom {{{\bf{F}}\left( {{\bf{x}},t} \right)\Delta t} {\rho \left( {{\bf{x}},t} \right)}}} \right.
\kern-\nulldelimiterspace} {\rho \left( {{\bf{x}},t} \right)}}$, and ${\bf{F}}\left( {{\bf{x}},t} \right)$ is the total force acting on a fluid node which is defined as ${\bf{F}}\left( {{\bf{x}},t} \right) = {{\bf{F}}_{{\mathop{\rm int}} }}\left( {{\bf{x}},t} \right) + {{\bf{F}}_{{\rm{ext}}}}\left( {{\bf{x}},t} \right)$ with ${{\bf{F}}_{{\mathop{\rm int}} }}$ and ${{\bf{F}}_{{\mathop{\rm ext}} }}$ representing the intermolecular force and the external force, respectively. In this setting, the real fluid velocity is given by
\begin{equation}
\rho {\bf{u}} = \sum\limits_{i = 0}^{18} {{{\bf{c}}_i}{f_i}}  + \frac{1}{2}{\bf{F}}\Delta t.
\label{eq5}
	\end{equation}
In this work, we assume that intermolecular forces act between pairs of molecules and are additive, then the interaction force can be given by \cite{ShanPRE2008}
\begin{equation}
{{\bf{F}}_{{\mathop{\rm int}} }} =  - G\psi \left( {\bf{x}} \right)\sum\limits_{i = 1}^{18} {\varpi \left( {{{\left| {{{\bf{c}}_i}} \right|}^2}} \right)} \psi \left( {{\bf{x}} + {{\bf{c}}_i}\Delta t} \right){{\bf{c}}_i},
\end{equation}
where $G$ is a simple scalar that controls the strength of the interaction, and it is usually set to -1, ${\varpi \left( {{{\left| {{{\bf{c}}_i}} \right|}^2}} \right)}$ is the weight given by $\varpi \left( 1 \right) = {1 \mathord{\left/{\vphantom {1 6}} \right. \kern-\nulldelimiterspace} 6}$, $\varpi \left( 2 \right) = {1 \mathord{\left/{\vphantom {1 {12}}} \right. \kern-\nulldelimiterspace} {12}}$. $\psi \left( {\bf{x}} \right)$ represents an effective mass (rather than the fluid density $\rho$) depending on the local density, and it is also named as the pseudopotential \cite{ShanPRE1993}. Following the work of He and Doolen \cite{HeJSP2002}, the form of the pseudopotential is suggested as 
\begin{equation}
\psi  = \sqrt {\frac{{2\Delta {t^2}\left( {{p_{EOS}} - \rho \hat c_s^2} \right)}}{{G\Delta {x^2}}}} .
\end{equation}
where $p_{EOS}$ is a prescribed non-ideal equation of state in the thermodynamic theory.      

Note that the evolution equation in LB community usually consists of two parts: collision and streaming. However, unlike the BGK model, the collision step in MRT model is executed in the moment space \cite{LiCF2019}, 
\begin{equation}
{{{\bf{\hat m}}}^*}\left( {{\bf{x}},t} \right) = {\bf{\hat m}}\left( {{\bf{x}},t} \right) - {\bf{\hat \Lambda }}\left[ {{\bf{\hat m}}\left( {{\bf{x}},t} \right) - {{{\bf{\hat m}}}^{eq}}\left( {{\bf{x}},t} \right)} \right] + \left[ {{{{\bf{\hat m}}}^{eq}}\left( {{\bf{x}},t} \right)\left| {_{{{\bf{u}}^{eq}} + \Delta {\bf{u}}}} \right. - {{{\bf{\hat m}}}^{eq}}\left( {{\bf{x}},t} \right)\left| {_{{{\bf{u}}^{eq}}}} \right.} \right],
\end{equation}
whereas the streaming step is implemented in the velocity space,
\begin{equation}
{f_i}\left( {{\bf{x}} + {{\bf{c}}_i}\Delta t,t + \Delta t} \right) = {f_i}^*\left( {{\bf{x}},t} \right).
\end{equation} 
Here ${\bf{\hat m}}\left( {{\bf{x}},t} \right) = {\left[ {{{\hat m}_0}\left( {{\bf{x}},t} \right),{{\hat m}_1}\left( {{\bf{x}},t} \right), \ldots ,{{\hat m}_{18}}\left( {{\bf{x}},t} \right)} \right]^{\rm T}} = {\bf{\hat M}}\left[ {{f_0}\left( {{\bf{x}},t} \right),{f_1}\left( {{\bf{x}},t} \right), \ldots ,{f_{18}}\left( {{\bf{x}},t} \right)} \right]$ is the post-collision moment,
${f_i}^*\left( {{\bf{x}},t} \right)$ represents the distribution function after collision which can be obtained via the inverse matrix of ${{\bf{\hat M}}}$ as ${f_i}^*\left( {{\bf{x}},t} \right) = {{{\bf{\hat M}}}^{ - 1}}{{{\bf{\hat m}}}^*}\left( {{\bf{x}},t} \right)$, ${{{\bf{\hat m}}}^{eq}}\left( {{\bf{x}},t} \right) = {\left[ {\hat m_0^{\left( {eq} \right)}\left( {{\bf{x}},t} \right), \ldots ,\hat m_{18}^{\left( {eq} \right)}\left( {{\bf{x}},t} \right)} \right]^{\rm T}} = {\bf{\hat M}}\left[ {f_0^{\left( {eq} \right)}\left( {{\bf{x}},t} \right),, \ldots ,f_{18}^{\left( {eq} \right)}\left( {{\bf{x}},t} \right)} \right]^{\rm T}$is the equilibrium moment function given by 
\begin{equation}
\left\{ \begin{array}{l}
		{\bf{\hat m}}_{0 - 6}^{eq}\left( {{\bf{x}},t} \right) = \left[ {\rho ,\rho {u_x},\rho {u_y},\rho {u_z},\rho  + \rho {{\left| {\bf{u}} \right|}^2},\left. {\rho \left( {2u_x^2 - u_y^2 - u_z^2} \right),\rho \left( {u_y^2 - u_z^2} \right)} \right]^{\rm T}} \right.\\
		{\bf{\hat m}}_{7 - 15}^{eq}\left( {{\bf{x}},t} \right) = \left[ {\rho {u_x}{u_y},\rho {u_x}{u_z},\rho {u_y}{u_z},\rho \hat c_s^2{u_y},\rho \hat c_s^2{u_x},\rho \hat c_s^2{u_z},\rho \hat c_s^2{u_x},\rho \hat c_s^2{u_z},\rho \hat c_s^2{u_y}} \right]^{\rm T}\\
		{\bf{\hat m}}_{16 - 18}^{eq}\left( {{\bf{x}},t} \right) = \left[ {\rho \hat c_s^4\left( {1 - 1.5{{\left| {\bf{u}} \right|}^2}} \right) + \rho \hat c_s^2\left( {u_x^2 + u_y^2} \right),\rho \hat c_s^4\left( {1 - 1.5{{\left| {\bf{u}} \right|}^2}} \right) + \rho \hat c_s^2\left( {u_x^2 + u_z^2} \right),\rho \hat c_s^4\left( {1 - 1.5{{\left| {\bf{u}} \right|}^2}} \right) + \rho \hat c_s^2\left( {u_y^2 + u_z^2} \right)} \right]^{\rm T}
	\end{array} \right..
\end{equation}
${{\bf{\hat \Lambda }}}$ is the diagonal relaxation matrix defined as
\begin{equation}
{\bf{\hat \Lambda }} = diag\left[ {1,1,1,1,{{\hat s}_e},{{\hat s}_\upsilon },{{\hat s}_\upsilon },{{\hat s}_\upsilon },{{\hat s}_\upsilon },{{\hat s}_\upsilon },{{\hat s}_q},{{\hat s}_q},{{\hat s}_q},{{\hat s}_q},{{\hat s}_q},{{\hat s}_q},{{\hat s}_\pi },{{\hat s}_\pi },{{\hat s}_\pi }} \right],
\end{equation}
where ${{{\hat s}_e}}$ and ${{{\hat s}_\upsilon }}$ rely on the bulk viscosity and shear viscosity, whereas ${{{\hat s}_q}}$ and ${{{\hat s}_\pi }}$ depend on the non-hydrodynamic moments.  Based on the above analysis, it is clear that the evolution process for the MRT model includes the transformation between the moment space and the velocity space, which must be careful treatment in practical programming.

\subsection{Lattice Boltzmann equation for temperature function} 
By neglecting the viscous heat dissipation, the temperature equation for non-ideal gases can be derived from the local balance law as \cite{ZhangPRE2021,LiPRE2022}
\begin{equation}
	\rho {c_v}\frac{{\partial T}}{{\partial t}} + \rho {c_v}{\bf{u}} \cdot \nabla T = \nabla  \cdot \left( {\lambda \nabla T} \right) - T{\left( {\frac{{\partial {p_{EOS}}}}{{\partial T}}} \right)_\rho }\nabla  \cdot {\bf{u}}.
	\label{eq12}
\end{equation}
where $T$, $\lambda$, and $c_v$ are the temperature, thermal conductivity, and specific heat at constant volume, respectively. To treat the liquid-vapor phase change problem, the above governing equation is usually rewritten as a standard convection-diffusion equation in previous LB models \cite{ZhangPRE2021,LiPRE2017,LiPRE2022}, such that a discrete source term related to $\nabla \left( {\rho {c_v}} \right)$ was also needed to be added into the temperature evolution, which may affect the numerical stability and accuracy when the fluid properties change significantly across the liquid-vapor interface \cite{HuATE2019,HuCF2019}. In this setting, recently we proposed a 2D MRT thermal LB model for liquid-vapor phase change \cite{WangPRE2022}, and in contrast to previous LB models, the calculation of $\nabla \left( {\rho {c_v}} \right)$ is eliminated in the proposed model. In this work,  considering the significant differences in the development and implementation of 3D MRT model, we extend this 2D model to 3D situation, and the evolution equation for the temperature $T$ distribution function ${g_i}\left( {{\bf{x}},t} \right)$ can be expressed as 
\begin{equation}
	\rho c_v{g_i}\left( {{\bf{x}} + {{\bf{c}}_i}\Delta t,t + \Delta t} \right) =
	{g_i}\left( {{\bf{x}},t} \right)  +(\rho c_v-1){g_i}\left( {{\bf{x}} + {{\bf{c}}_i}\Delta t,t} \right)
	-\left( {{{{\bf{ M}}}^{ - 1}}\Lambda {\bf{ M}}} \right)_{ij}\left[
	{{g_j}\left( {{\bf{x}},t} \right) - g_j^{(eq)}\left( {{\bf{x}},t}
		\right)} \right]+\Delta t {\bar F_i} + \vartheta \Delta t {S_i}  ,
	\label{eq13}
\end{equation} 
where ${\bf{ M}}$ is the non-orthogonal transformation matrix given in Appendix \cite{HeIJHMT2019,LiCF2019}, $g_i^{\left( {eq} \right)}$ is the temperature equilibrium distribution and it is a linear function 
\begin{equation}
	g_i^{(eq)} = {{w}_i}T,
	\label{eq14}
\end{equation} 
where ${w_0} = 1 - \overline w ,{w_{1 - 6}} = {{\overline w } \mathord{\left/
{\vphantom {{\overline w } 6}} \right.	\kern-\nulldelimiterspace} 6}$ with $\overline w  \in \left( {0,1} \right)$, and the speed of sound for D3Q7 lattice is given by $c_s^2 = {{\overline w } \mathord{\left/	{\vphantom {{\overline w } 3}} \right.\kern-\nulldelimiterspace} 3}$ \cite{LiPRE2022}. $G_i$ is the source term defined as 
\begin{equation}
	{\bar F_i} =- {{ w}_i}{\bar F},
	\label{eq15}
\end{equation}
where $\bar F$ is the forcing term given by ${\bar F} =  {\rho {c_v}{\bf{u}} \cdot \nabla T + T{{\left( {{\partial _T}{p_{EOS}}} \right)}_\rho }\nabla  \cdot {\bf{u}}} $. $S_i$ is a small correction term aimed at eliminating the additional terms in the recovered temperature equation   
\begin{equation}
	S_i={{w}_i} \rho c_v\frac{\Delta t}{2} \partial_t^2T,
	\label{eq16}
\end{equation}  
and $\vartheta $ is a free parameter used to reflect the influence of the correction term, and its value may be 0 or 1, which will be discussed later. $\Lambda$ is the diagonal relaxation matrix expressed as  $\Lambda  = \left( {{\varsigma _0},{\varsigma _1}, \ldots ,{\varsigma _6}} \right)$ with ${\varsigma _i} \in \left( {0,2} \right)$. 

We now turn to carry out the Chapman-Enskog analysis for the present model to show the temperature equation can be recovered. Applying the Taylor series expansion to Eq. (\ref{eq13}), we have     
\begin{equation}
	\rho c_v\left[g_i+ \Delta t D_ig_i+\frac{\Delta t^2}{2}D_i^2g_i    \right] -
	g_i  +(1-\rho c_v)\left[g_i+\Delta td_i g_i +\frac{\Delta t^2}{2}d_i^2g_i  \right]
	=- {\left( {{{\bf{M}}^{ - 1}}\Lambda {\bf{M}}} \right)_{ij}} \left[{g_i - g_i^{(eq)}} \right]   +\Delta t {\bar F_i}+\vartheta \Delta t {S_i},
	\label{eq17}
\end{equation}
Here, ${D_i} = {\partial _t} + {d_i}$ with ${d_i} = {{\bf{c}}_i} \cdot \nabla$. Substituting the following Chapman-Enskog expansions 
\begin{equation}
{\partial _t} = \sum\limits_{n = 1}^{ + \infty } {{\varepsilon ^n}\partial {t_n}} , \; \;\;\nabla  = \varepsilon {\nabla _1}, \; \;\;{g_i} = \sum\limits_{n = 1}^{ + \infty } {{\varepsilon ^n}g_i^{\left( n \right)}} , \; \;\; \bar F = \varepsilon {\bar F^{\left( 1 \right)}},
\end{equation}
into Eq. (\ref{eq17}), and separating the results into terms of different order, we can obtain 

\begin{equation}
	{\rm O}(\epsilon ^0): \quad g_i^{(0)}=g_i^{(eq)},
	\label{eq19}
\end{equation}
\begin{equation}
	{\rm O}(\epsilon ^1):\quad \rho c_v\partial_ {t_1}g_i^{(0)}+d_{1i}g_i^{(0)}=-\frac{1}{\Delta t} \left({{{\bf{ M}}}}^{-1}{\Lambda}{{{\bf{ M}}}}\right)_{ij} g_j^{(1)} + {{w}_i \bar F^{(1)}},
	\label{eq42}
\end{equation}
\begin{equation}
	{\rm O}(\epsilon ^2): \quad \rho c_vD_{1i}g_i^{(1)}+\rho c_v\partial_ {t_2}g_i^{(0)}+\rho c_v \frac{\Delta t}{2}D_{1i}^2g_i^{(0)}+(1-\rho c_v)d_{1i}g_i^{(1)}+(1-\rho c_v)\frac{\Delta t}{2} d_{1i}^2g_i^{(0)}=-\frac{1}{\Delta t} \left({{{\bf{ M}}}}^{-1}{\Lambda}{{{\bf{ M}}}}\right)_{ij}g_j^{(2)}+{ w}_i \vartheta  \rho c_v\frac{\Delta t}{2} \left(\partial t_1\right)^2T.
	\label{eq43} 
\end{equation}
where ${D_{{1i}}} = {\partial _{{t_1}}} + {d_{{1i}}} = {\partial _{{t_1}}} + {{\bf{c}}_i} \cdot {\nabla _1}$. Then the corresponding 
equations in the moment space can be obtained by multiplying the non-orthogonal transformation matrix ${\bf{ M}}$ on both sides of the above three equations, 
\begin{equation}
	{\rm O}(\epsilon ^0):\quad {{{\bf{ m}}}}^{(0)}={{{\bf{ m}}}}^{(eq)},
	\label{eq22} 
\end{equation}
\begin{equation}
	{\rm O}(\epsilon ^1):\quad \rho c_v {\bf{I}} \partial_ {t_1}{{{\bf{ m}}}}^{(0)}+{\bf{d}}_{1}{{{\bf{ m}}}}^{(0)}=-\frac{1}{\Delta t} {\Lambda}{{{\bf{ m}}}} ^{(1)} + {{{\bf{ M}}}}{{{\bf{ \bar F}}}}^{(1)},
	\label{eq23} 
\end{equation}
\begin{equation}
	O(\epsilon ^2): \quad \rho c_v{\bf{I}}\partial_{t_1}  {{{\bf{ m}}}}^{(1)}+  {\bf{d}}_{1}{{{\bf{ m}}}}^{(1)} +\rho c_v{\bf{I}}\partial_ {t_2} {{{\bf{ m}}}}^{(0)}+\rho c_v{\bf{I}} \frac{\Delta t}{2}\partial_{t_1}^2{{{\bf{ m}}}}^{(0)}+\rho c_v  \frac{\Delta t}{2}\left(2\partial_{t_1}{\bf{d}}_{1}\right){{{\bf{ m}}}}^{(0)}+ \frac{\Delta t}{2}{\bf{d}}_{1}^2{{{\bf{ m}}}}^{(0)}=- \frac{1}{\Delta t} {\Lambda } {{{\bf{m}}}}^{(2)}+ \vartheta {\bf{M}}{{{\bf{ S}}}}^{(2)},
	\label{eq24} 
\end{equation}
where ${{\bf{d}}_1} = {\bf{ M}}{\rm{diag}}\left[ {{{\bf{c}}_0} \cdot \nabla , \ldots ,{{\bf{c}}_6} \cdot \nabla } \right]{{{\bf{ M}}}^{ - 1}}$, ${{{\bf{ \bar F}}}^{\left( 1 \right)}} = {\left[ {{{ w}_0}{\bar F^{\left( 1 \right)}}, \ldots ,{{ w}_6}{\bar F^{\left( 1 \right)}}} \right]^{\rm T}}$, ${{{\bf{ S}}}^{\left( 2 \right)}} = {\left[ {0.5{{ w}_0}\rho {c_v}\partial _{{t_1}}^2T, \ldots ,0.5{{ w}_6}\rho {c_v}\partial _{{t_1}}^2T} \right]^{\rm T}}$. ${\bf{ m}} = {\bf{ Mg}}$ is the moment of temperature distribution function, and ${{{\bf{ m}}}^{\left( {eq} \right)}}$ in Eq. (\ref{eq22}) is the equilibrium moment function defined as
\begin{equation}
	{{{\bf{ m}}}^{\left( {eq} \right)}} = {\bf{M}}{{\bf{g}}^{\left( {eq} \right)}} = {\left[ {T,0,0,0,\overline w T,0,0} \right]^{\rm T}}.
	\label{eq47} 
\end{equation} 
According to Eq. (\ref{eq23}), we can obtain the following first-order equations of the  conserved moments ($m_0$, $m_1$, $m_2$ and $m_3$), 
\begin{equation}	
	\rho c_v \partial_ {t_1}T= {\bar F}^{(1)},
	\label{eq26}
\end{equation}
\begin{equation}
	c_s^2\partial _{x1}T =-\frac{1}{\Delta t} {\varsigma}_1{m}_1 ^{(1)},
	\label{eq27}
\end{equation}
\begin{equation}
	c_s^2\partial _{y1}T =-\frac{1}{\Delta t} {\varsigma}_2{m}_2 ^{(1)},
	\label{eq28}
\end{equation}
\begin{equation}
	c_s^2\partial _{z1}T =-\frac{1}{\Delta t} {\varsigma}_3{m}_3 ^{(1)},
	\label{eq29}
\end{equation}
  Note that $\mathbf{d}_{1}^{2} \mathbf{m}^{(0)}=\left(\mathbf{E} \cdot \nabla_{1}\right)\left(\mathbf{E} \cdot \nabla_{1}\right) \mathbf{m}^{(0)}$ where $\mathbf{E}=\left(\mathbf{E}_{x}, \mathbf{E}_{y}, \mathbf{E}_{z}\right)$, and the details of the elements in ${\bf{E}}$  are given as

\begin{equation}
\left\{ \begin{array}{l}
	{{\bf{E}}_x} = {\bf{M}}{\rm{diag}}\left[ {{{\bf{c}}_{0x}}, \ldots ,{{\bf{c}}_{6x}}} \right]{{\bf{M}}^{ - 1}}\\
	{{\bf{E}}_y} = {\bf{M}}{\rm{diag}}\left[ {{{\bf{c}}_{0y}}, \ldots ,{{\bf{c}}_{6y}}} \right]{{\bf{M}}^{ - 1}}\\
	{{\bf{E}}_z} = {\bf{M}}{\rm{diag}}\left[ {{{\bf{c}}_{0z}}, \ldots ,{{\bf{c}}_{6z}}} \right]{{\bf{M}}^{ - 1}}
\end{array} \right.,
\end{equation}
then we can  obtain the following second-order equation,
\begin{equation}
	\rho {c_v}{\partial _{{t_2}}}T + {\partial _{x1}}\left\{ {{{ m}^{\left( 1 \right)}} + \frac{{\Delta t}}{2} c_s^2{\partial _{x1}}T} \right\} + {\partial _{y1}}\left\{ {{{ m}^{\left( 2 \right)}} + \frac{{\Delta t}}{2} c_s^2{\partial _{y1}}T} \right\} +  {\partial _{z1}}\left\{ {{{ m}^{\left( 3 \right)}} + \frac{{\Delta t}}{2} c_s^2{\partial _{z1}}T} \right\} +\left( {1 - \vartheta } \right)\frac{{\Delta t}}{2}\rho {c_v}\partial _{{t_1}}^2T = 0.
	\label{eq511}
\end{equation}
Based on the equations of Eqs. (\ref{eq26})-(\ref{eq29}), we have  
\begin{equation}
	\rho {c_v}{\partial _{{t_2}}}T = {\partial _{x1}}\left\{ {\Delta t\left( {\frac{1}{{{\varsigma _1}}} - \frac{1}{2}} \right) c_s^2{\partial _{x1}}T} \right\} + {\partial _{y1}}\left\{ {\Delta t\left( {\frac{1}{{{\varsigma _2}}} - \frac{1}{2}} \right) c_s^2{\partial _{y1}}T} \right\}  + {\partial _{z1}}\left\{ {\Delta t\left( {\frac{1}{{{\varsigma _3}}} - \frac{1}{2}} \right) c_s^2{\partial _{z1}}T} \right\} + \left( {1 - \vartheta } \right)\frac{{\Delta t}}{2}\rho {c_v}\partial _{{t_1}}^2T,
	\label{eq32}
\end{equation} 
Considering Eq. (\ref{eq26}) and Eq. (\ref{eq32}), the temperature equation is recovered as   
\begin{equation}
	\rho {c_v}{\partial _{{t}}}T = \nabla  \cdot \left( {\lambda \nabla T} \right) + {\bar F} + \left( { \vartheta -1 } \right)\frac{{\Delta t}}{2}\rho {c_v}\partial _{{t}}^2T  ,	
\end{equation} 
where $\lambda  = \Delta t\left( {{\varsigma _1}^{ - 1} - 0.5} \right) c_s^2 = \Delta t\left( {{\varsigma _2}^{ - 1} - 0.5} \right)c_s^2 = \Delta t\left( {{\varsigma _3}^{ - 1} - 0.5} \right)c_s^2$.

The above multi-scale analysis shows that the temperature equation can be correctly recovered at the order of ${\varepsilon ^2}$ if the free parameter $\vartheta$ is set to 1. However, this treatment requires an additional finite-difference scheme to calculate $\partial _t^2T$, and this, in turn, increases the computational resource, which is undesirable for 3D numerical simulations. Fortunately, as shown below, the numerical results depend little on the correction term $S_i$. Thus, to make the present model more efficient, the value of the free parameter $\vartheta$ is recommended to be zero in practical applications. In such case, the Chapman-Enskog analysis is equivalent to performing at second-order in space and first-order in time. Moreover, the multi-scale analysis shows that the temperature gradient $\nabla T$ appeared in the discrete forcing term (see Eq. (\ref{eq15})) can be calculated by a local scheme shown in Eqs. (\ref{eq27})-(\ref{eq29}), while for the gradient term $\nabla  \cdot {\bf{u}}$ in ${\overline F _i}$ , it is calculated by using an isotropic finite-difference scheme \cite{LeeJCP2005}. 

\section{Boundary condition treatment}
Unlike the traditional numerical method in which the boundary condition is specified in terms of the fluid variables such as velocity, pressure, and temperature, the primitive variable in the LB method is the discrete distribution. Therefore, it is important for the LB method to find a way to transform the boundary condition from the fluid variables to those for distribution functions. In this section,  we will present some necessary boundary condition treatments for simulating 3D liquid-vapor phase change.                  
\subsection{Macroscopic variables}
For the velocity field, the modified bounce-back scheme is adopted for the stationary nonslip wall. In this scheme, the collision process is  imposed on the boundary, and the unknown pre-collision distributions at a boundary node ${\bf x}_b$ is expressed as \cite{LaddJFM1994} 
\begin{equation}
{f_{\bar i}}\left( {{{\bf{x}}_b},t + \Delta t} \right) = f_i^*\left( {{{\bf{x}}_b},t} \right),
\end{equation}
where ${f_{\bar i}}$ refers to the distribution function for the direction $\bar i$ such that ${{\bf{c}}_{\bar i}} =  - {{\bf{c}}_i}$. After this treatment, all the distribution functions at the boundary nodes are obtained, and then the corresponding fluid density can be calculated from $\rho \left( {{{\bf{x}}_b},t} \right) = \sum\limits_{i = 0}^{18} {{f_i}\left( {{{\bf{x}}_b},t} \right)} $, while the velocity ${\bf{u}}\left( {{{\bf{x}}_b},t} \right)$ is directly specified according to the nonslip boundary condition. Obviously, the aforementioned scheme is capable of exactly satisfying mass conservation at a boundary node. 

As for the temperature field,  here we only consider the Dirichlet boundary condition. To obtain the unknown distributions at the boundary node ${{{\bf{x}}_b}}$, the non-equilibrium extrapolation scheme \cite{GuoPOF2002} is adopted. The basic idea of this scheme is to decompose the distribution function at a boundary node into its equilibrium and non-equilibrium parts, 
\begin{equation}
{g_i}\left( {{{\bf{x}}_b},t} \right) = g_i^{\left( {eq} \right)}\left( {{{\bf{x}}_b},t} \right) + g_i^{\left( {neq} \right)}\left( {{{\bf{x}}_b},t} \right),
\end{equation}  
where $g_i^{\left( {eq} \right)}\left( {{{\bf{x}}_b},t} \right)$ is directly  obtained according to the temperature given by the  Dirichlet boundary condition, while $g_i^{\left( {neq} \right)}\left( {{{\bf{x}}_b},t} \right)$ is approximated by the non-equilibrium part at the neighboring fluid node ${{\bf{x}}_f}$,
\begin{equation}
g_i^{\left( {neq} \right)}\left( {{{\bf{x}}_b},t} \right) = {g_i}\left( {{{\bf{x}}_f},t} \right) - g_i^{\left( {eq} \right)}\left( {{{\bf{x}}_f},t} \right).
\end{equation} 
It has been demonstrated that the above extrapolation scheme is of second-order accuracy in both time and space. 

\subsection{Wetting boundary condition}

\begin{figure}[H]
	\centering
		\includegraphics[width=0.45\textwidth]{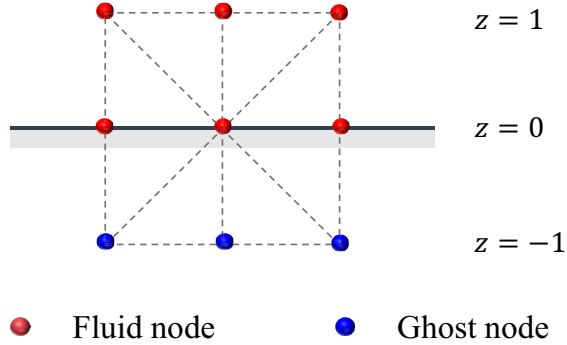}\\		
	\caption{Sketch of the geometric formulation scheme. }
	\label{fig2}
\end{figure}

When modeling multiphase flow, apart from the common boundary conditions for the velocity and temperature fields, we also require an additional wetting boundary condition. In this work, the pseudopotential-based geometric formulation is used to implement the wetting boundary condition for its better performance in simulations of super-hydrophilic or super-hydrophobic dynamics. As sketched in Fig. 2,  a ghost layer is put outside the wall in this approach, and the pseudopotential $\psi \left( {\bf{x}} \right)$ at this layer is given by \cite{WuPRE2020}
\begin{equation}
{\psi _{x,y,-1}} = {\psi _{x,y,1}} + \sqrt {{{\left( {{\psi _{x + 1,y,0}} - {\psi _{x - 1,y,0}}} \right)}^2} + {{\left( {{\psi _{x,y + 1,0}} - {\psi _{x,y - 1,0}}} \right)}^2}} \tan \left( {\frac{\pi }{2} - \theta } \right),
\end{equation}   
where $\theta $ is the prescribed contact angle. 
\subsection{Outflow boundary condition}
\label{sec41}
Since the physical domain that appeared in practical multiphase flow problems is usually very large, it is necessary to introduce proper outflow boundary conditions at the open boundary, and thus it allows us to use a reduced computational domain. As for the liquid-vapor phase change, the most commonly used outflow boundary condition is the so-called convective boundary condition (CBC) \cite{LouPRE2013}; however, our numerical tests show that the choice of the typical velocity significantly influences the numerical performance. In this setting, inspired by the Zhou-He Scheme \cite{ZouPOF1997,WangIJNMF2020}, in what follows we intend to develop a constant-pressure boundary condition on the D3Q19 lattice to model the outflow boundary condition,  and we note it has not yet been reported in previous works.     

Let us assume that the flow is in the z-direction, then the unknown density distribution functions at the outlet are ${f_6},{f_{12}},{f_{13}},{f_{16}}$ and ${f_{17}}$, which need to be determined. According to Eq. (\ref{eq5}), we have   
\begin{equation}
		\rho {u_x} = \sum\limits_{i = 0}^{18} {{c_{ix}}{f_i} + \frac{{\Delta t}}{2}{F_x}}, 	
		 \label{eq38}
\end{equation}

\begin{equation}
		\rho {u_y} = \sum\limits_{i = 0}^{18} {{c_{iy}}{f_i} + \frac{{\Delta t}}{2}{F_y}}, 
 \label{eq39}
\end{equation}

\begin{equation}
	\rho {u_z} = \sum\limits_{i = 0}^{18} {{c_{iz}}{f_i} + \frac{{\Delta t}}{2}{F_z}}. 
	 \label{eq40}
\end{equation}

Suppose ${u_x}$ and ${u_y}$ at the outlet are both equal to zero, and note that the density $\rho$ at the outlet is known, we get 
\begin{equation}
{u_z} = \frac{1}{\rho }\left[ {{f_0} + {f_1} + {f_2} + {f_3} + {f_4} + {f_7} + {f_8} + {f_9} + {f_{10}} + 2\left( {{f_5} + {f_{11}} + {f_{14}} + {f_{15}} + {f_{18}}} \right) + \frac{1}{2}{F_z}} \right] - 1,
\end{equation}

Using the bounce back scheme for the non-equilibrium part of the distribution functions in the normal direction, i.e., ${f_5} - f_5^{\left( {eq} \right)} = {f_6} - f_6^{\left( {eq} \right)}$ we have 
\begin{equation}
	{f_6} = {f_5} - \frac{1}{3}\rho {u_z}.
	\label{eq42}
\end{equation}
To ensure the hold of Eq.(\ref{eq38}), we further modified the distribution functions of ${f_{12}},{f_{13}},{f_{16}}$ and ${f_{17}}$ as 
\begin{equation}
	{f_{12}} = {f_{11}} - f_{11}^{\left( {eq} \right)} + f_{12}^{\left( {eq} \right)} - \delta x - \delta z,
\end{equation} 
\begin{equation}
	{f_{13}} = {f_{14}} - f_{14}^{\left( {eq} \right)} + f_{13}^{\left( {eq} \right)} + \delta x - \delta z,
\end{equation} 
\begin{equation}
	{f_{16}} = {f_{15}} - f_{15}^{\left( {eq} \right)} + f_{16}^{\left( {eq} \right)} - \delta y - \delta z,
\end{equation} 
\begin{equation}
	{f_{17}} = {f_{18}} - f_{18}^{\left( {eq} \right)} + f_{17}^{\left( {eq} \right)} + \delta y - \delta z
\end{equation} 
Submitting the above equations into Eqs.(\ref{eq38})-(\ref{eq40}) yields
\begin{equation}
	\delta x =  - \frac{1}{2}\left( {{f_1} - {f_2}} \right) - \frac{1}{2}\left( {{f_7} - {f_8}} \right) - \frac{1}{2}\left( {{f_9} - {f_{10}}} \right) - \frac{1}{4}{F_x},
\end{equation}
\begin{equation}
		\delta y =  - \frac{1}{2}\left( {{f_3} - {f_4}} \right) - \frac{1}{2}\left( {{f_7} - {f_8}} \right) + \frac{1}{2}\left( {{f_9} - {f_{10}}} \right) - \frac{1}{4}{F_y},	
\end{equation}
\begin{equation}
	\delta z =  - \frac{1}{8}{F_z}.
\end{equation}
Then the distribution functions of  ${f_{12}},{f_{13}},{f_{16}}$ and ${f_{17}}$ can be expressed as 
\begin{equation}
	{f_{12}} = {f_{11}} + \frac{1}{2}\left( {{f_1} - {f_2}} \right) + \frac{1}{2}\left( {{f_7} - {f_8}} \right) + \frac{1}{2}\left( {{f_9} - {f_{10}}} \right) + \frac{1}{8}\left( {2{F_x} + {F_z}} \right) - \frac{1}{6}\rho {u_z},
\end{equation}
\begin{equation}
	{f_{13}} = {f_{14}} - \frac{1}{2}\left( {{f_1} - {f_2}} \right) - \frac{1}{2}\left( {{f_7} - {f_8}} \right) - \frac{1}{2}\left( {{f_9} - {f_{10}}} \right) - \frac{1}{8}\left( {2{F_x} - {F_z}} \right) - \frac{1}{6}\rho {u_z},
\end{equation}

\begin{equation}
	{f_{16}} = {f_{15}} + \frac{1}{2}\left( {{f_3} - {f_4}} \right) + \frac{1}{2}\left( {{f_7} - {f_8}} \right) - \frac{1}{2}\left( {{f_9} - {f_{10}}} \right) + \frac{1}{8}\left( {2{F_y} + {F_z}} \right) - \frac{1}{6}\rho {u_z},
\end{equation}

\begin{equation}
			{f_{17}} = {f_{18}} - \frac{1}{2}\left( {{f_3} - {f_4}} \right) - \frac{1}{2}\left( {{f_7} - {f_8}} \right) + \frac{1}{2}\left( {{f_9} - {f_{10}}} \right) - \frac{1}{8}\left( {2{F_y} - {F_z}} \right) - \frac{1}{6}\rho {u_z}.
\end{equation}
Together with Eq.(\ref{eq42}), all the unknown density distribution functions are established. 

\section{Validations and discussions}
In this section, we consider three numerical tests to validate the capacity of the proposed 3D thermal LB model in simulating liquid-vapor phase change. The Peng-Robinson (PR) equation of state is used in our simulations, and it is given by \cite{YuanPOF2006}
\begin{equation}
	{p_{EOS}} = \frac{{\rho RT}}{{1 - b\rho }} - \frac{{a\xi \left( T \right){\rho ^2}}}{{1 + 2b\rho  - {b^2}{\rho ^2}}},
\end{equation}  
where $\xi \left( T \right) = {\left[ {1 + \left( {0.37464 + 1.54226\omega  - 0.26992{\omega ^2}} \right)\left( {1 - \sqrt {{T \mathord{\left/{\vphantom {T {{T_c}}}} \right.\kern-\nulldelimiterspace} {{T_c}}}} } \right)} \right]^2}$ with $\omega$ and $T_c$ being the acentric factor and the critical temperature, respectively. $a$ and $b$ are two built-in variables defined as  
\begin{equation}
	a = \frac{{0.45724{R^2}T_c^2}}{{{p_c}}},b = \frac{{0.0778R{T_c}}}{{{p_c}}},
\end{equation}
where $R$ is the gas constant, $p_c$ is the critical pressure. In the following simulations, the basic parameters in PR equation of state are set as $a = {3 \mathord{\left/{\vphantom {3 {49}}} \right.\kern-\nulldelimiterspace} {49}},b = {2 \mathord{\left/{\vphantom {2 {21}}} \right.\kern-\nulldelimiterspace} {21}},R = 1,\omega  = 0.344$ \cite{LiPRE2017}. Additionally, the relaxation time ${{\hat s}_v}$ in the velocity field is equal to 1.25 (i.e., the shear viscosities for liquid and vapor phases are both fixed at 0.25), while the relaxation time  $\hat s_e $ is selected as 0.8. As for the relaxation time in the temperature field, the relaxation times ${\varsigma _1}$, ${\varsigma _2}$ and ${\varsigma _3}$ are determined by the conductivity. Apart from the above mentioned relaxation times, the rest relaxation times in velocity and temperature fields are taken as 1.
\subsection{Heat transfer in a 3D saturate liquid-vapor system}
\label{sec4120}
\begin{figure}[H]
	\centering
	\includegraphics[width=0.5\textwidth]{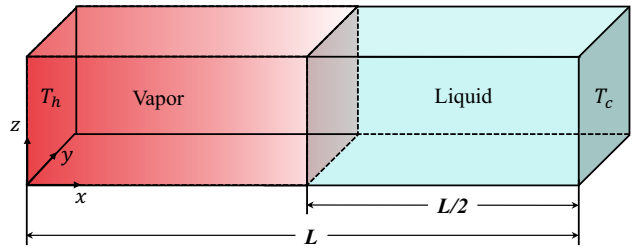}\\		
	\caption{Sketch of heat transfer in a 3D saturate liquid-vapor system. }
	\label{fig3}
\end{figure}
We first validate the present thermal model by considering the heat transfer in a 3D saturate liquid-vapor system \cite{GuoIJHMT2019}. As sketched in Fig. \ref{fig3}, initially  the vapor phase and liquid phase are uniformly distribute for $0 \le x \le {L \mathord{\left/ {\vphantom {L 2}} \right. \kern-\nulldelimiterspace} 2}$ and ${L \mathord{\left/{\vphantom {L 2}} \right.\kern-\nulldelimiterspace} 2} \le x \le L$ with the saturate temperature of $T = 0.86{T_c}$. At time $t>0$, a constant temperature $T_h=1.0T_c$ is imposed on the left wall, while the right wall is fixed at the saturate temperature of $T_{sat}$.  According to the Fourier's law and the continuity of heat flux, the slope ratio of the temperature profile should be equal to that of the thermal conductivities, and this phenomenon is independent of the choice of the specific heat $c_v$.

\begin{figure}[H]
	\centering
	\includegraphics[width=0.5\textwidth]{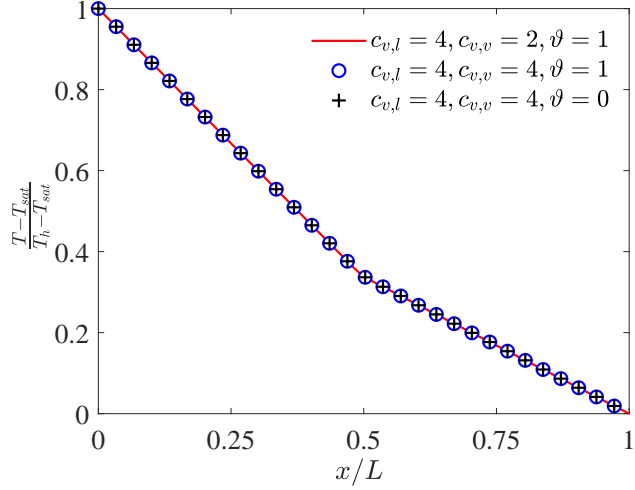}\\		
	\caption{Profiles of the temperature along the horizontal midplane ($x$-mid) for heat transfer in a 3D saturate liquid-vapor system. }
	\label{fig4}
\end{figure}

In simulations, the grid resolution for this problem is set to be $200 \times 10 \times 10$, and the periodic boundary condition is applied to the $y$ and $z$ directions. The thermal conductivity for liquid and vapor phases is taken as 0.4 and 0.2, respectively. Fig. \ref{fig4} presents the temperature distributions in the horizontal direction, where the effect of specific heat $c_v$ is also incorporated. It is found that the  slope ratio gained from the present numerical results is about 1.9507, which is very close to the theoretic result 2.0 (the relative error is $2.46\%$). In addition, the simulation results show that the variation of the specific heat at constant volume makes no difference to the temperature distribution, which is consistent with our previous analysis. Moreover, since the numerical results obtained from the case with and without the correction term  $S_i$ are almost the same, setting the free parameter in the temperature evolution equation equal to 0 in the actual application is appropriate. 

\subsection{Droplet evaporation in  an open space}
In this section, we intend to validate the proposed thermal LB model by simulating the droplet evaporation in an open space. For this problem, it turns out that the square of the droplet diameter tends to evolve linearly in time (also called $d^2$ law) \cite{WangPRE2022}, 
\begin{equation}
{\left( {{{{D \mathord{\left/	{\vphantom {D D}} \right.\kern-\nulldelimiterspace} D}}_0}} \right)^2} = 1 - Kt,
\end{equation} 
where $D_0$ is the droplet diameter at time $t=0$, $K$ is the evaporation rate. Numerical investigations are conducted in a cubic cavity with the grid size being $100 \times 100 \times 100$, and the periodic conditions are imposed in $x$, $y$ and $z$ direction for the velocity field. Initially, there is a droplet with the radius of $R_0 = 25 \Delta x$ placed at the center of the computational domain ${\bf{x}}_0$, and the temperature of the droplet is taken as $T_{sat}=0.86T_c$, while the gas temperature around the droplet is given by $T_{g}= 1.0T_c$. The velocity in the cubic cavity is equal to zero at time $t=0$, and the density field is initial as 
\begin{equation}
	\rho \left( {\bf{x}} \right) = \frac{{{\rho _l} + {\rho _v}}}{2} - \frac{{{\rho _l} - {\rho _v}}}{2}\tanh \frac{{2\left( {\left| {{\bf{x}} - {{\bf{x}}_0}} \right| - {R_0}} \right)}}{W},
\end{equation} 
where $W$ is the thickness of the phase interface, and it is set as $5 \Delta x$ in this work, $\rho_l$ and $\rho _v$ are the thermodynamic coexistence liquid and gas densities calculated by the Maxwell construction at $T_{sat}$. In simulations, the thermal conductivity $\lambda$ and the specific heat at constant volume $c_v$ are taken as $1.0/3.0$ and $5.0$, respectively. In addition, the temperature at boundaries is fixed at $T_g$, and the free parameter $\vartheta$ is chosen as 0. All the other unmentioned parameters used here are the same as that used in Sec. \ref{sec4120}.

\begin{figure}
	\centering 
	\subfigure{ 
		\includegraphics[width=0.25\textwidth]{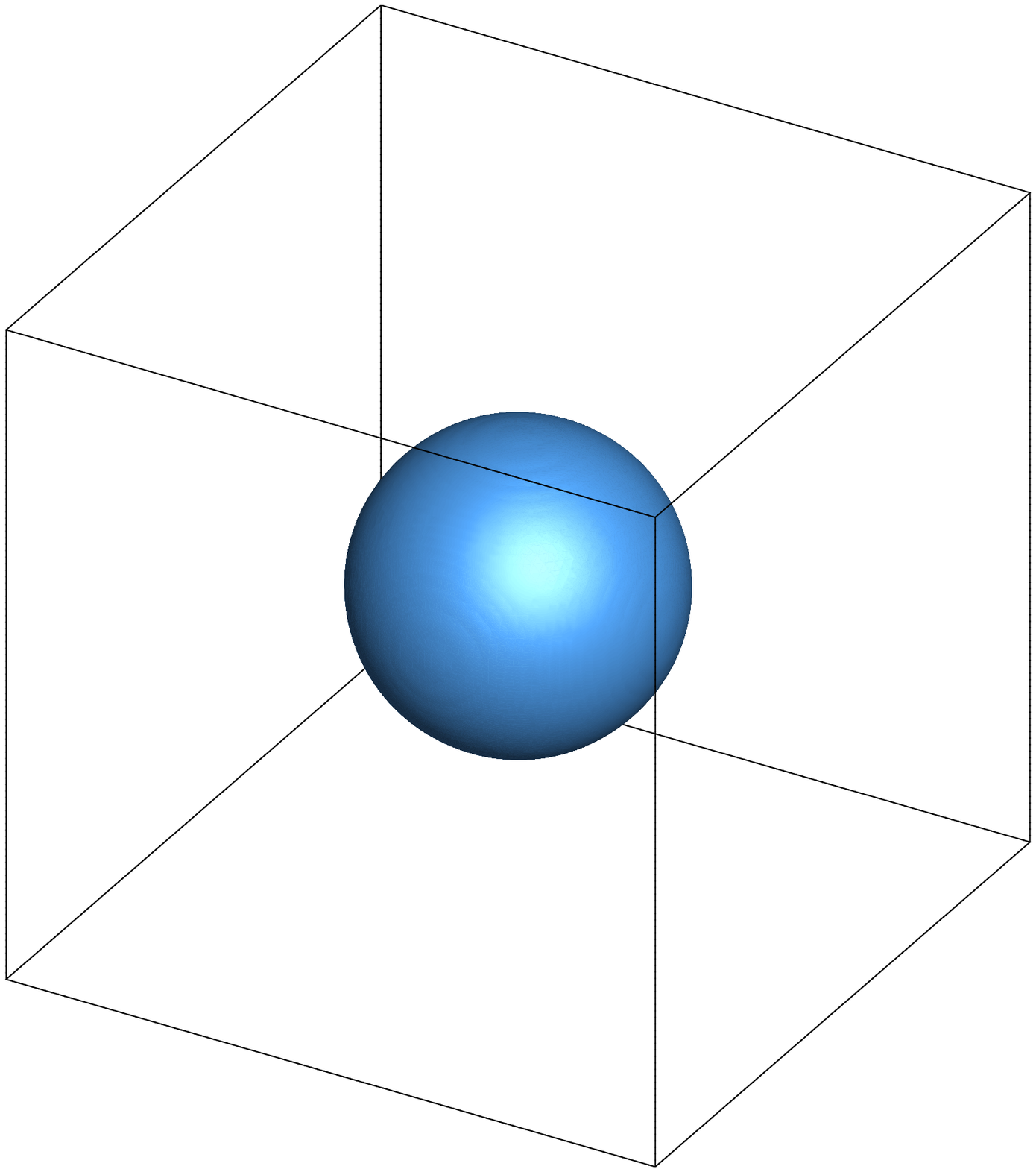}} 
	\subfigure{ 
		\includegraphics[width=0.25\textwidth]{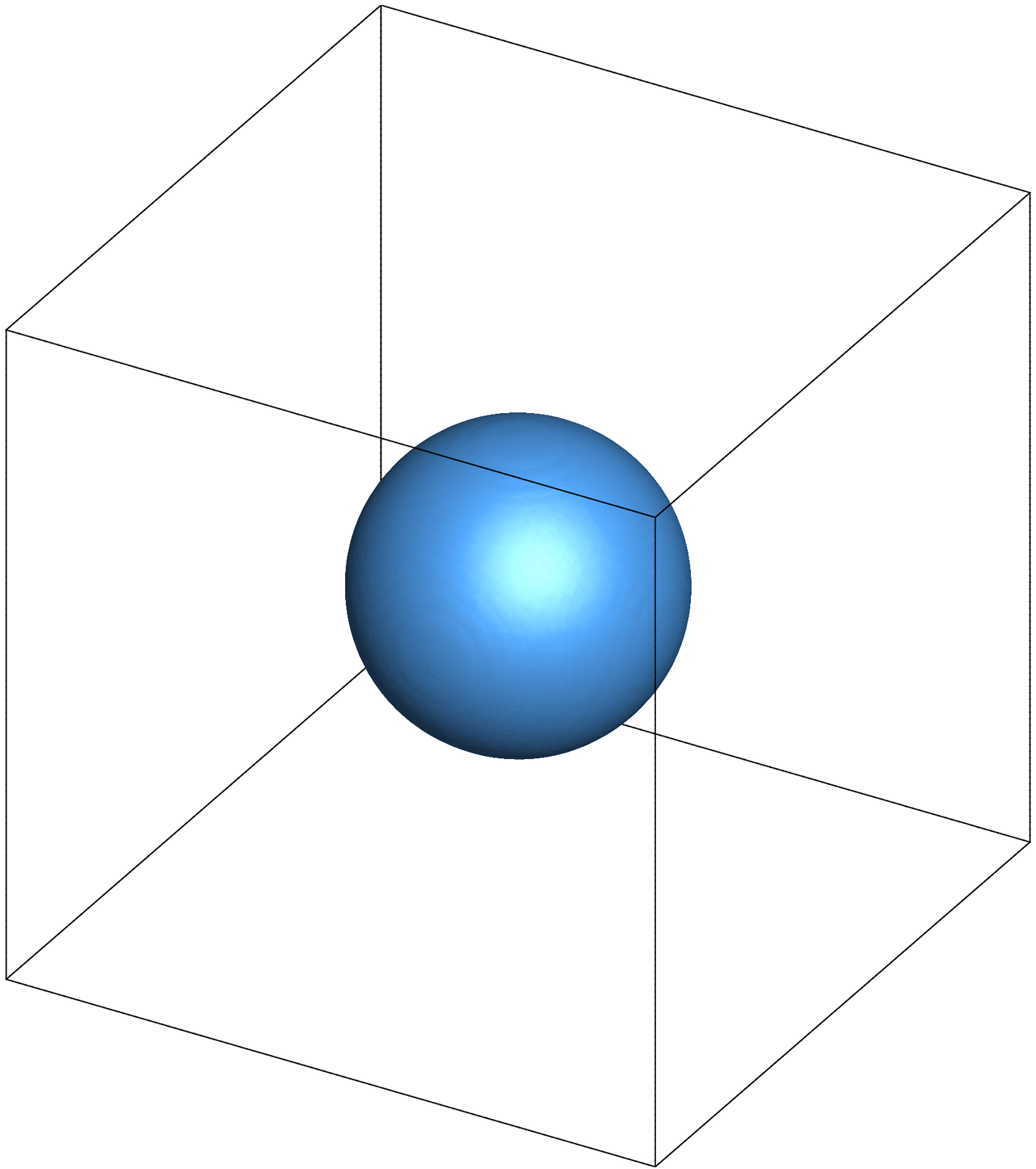}}
	\caption*{ (a) t = 40000$\delta t$}	
	\subfigure{ 
		\includegraphics[width=0.25\textwidth]{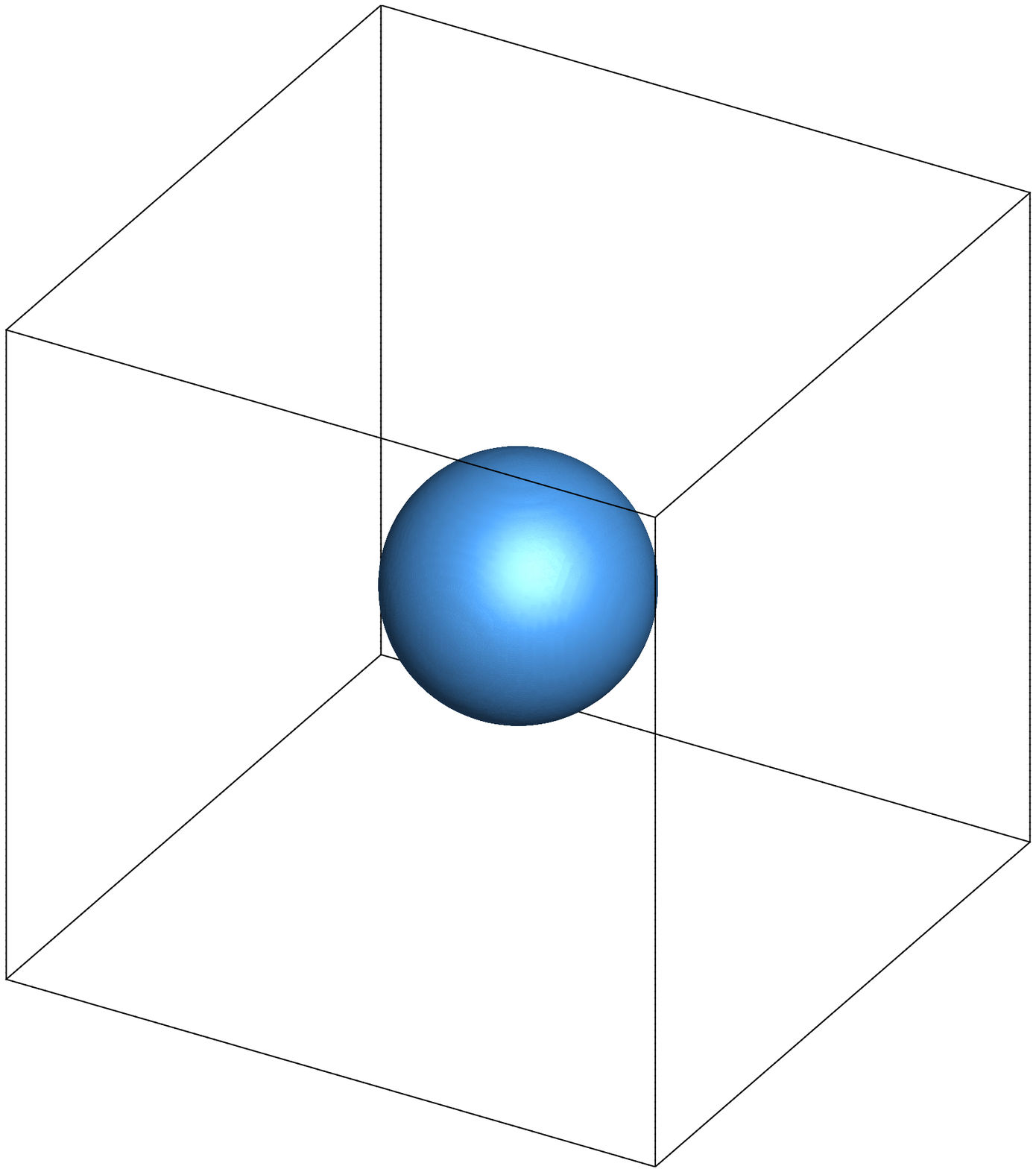}}	
	\subfigure{ 
		\includegraphics[width=0.25\textwidth]{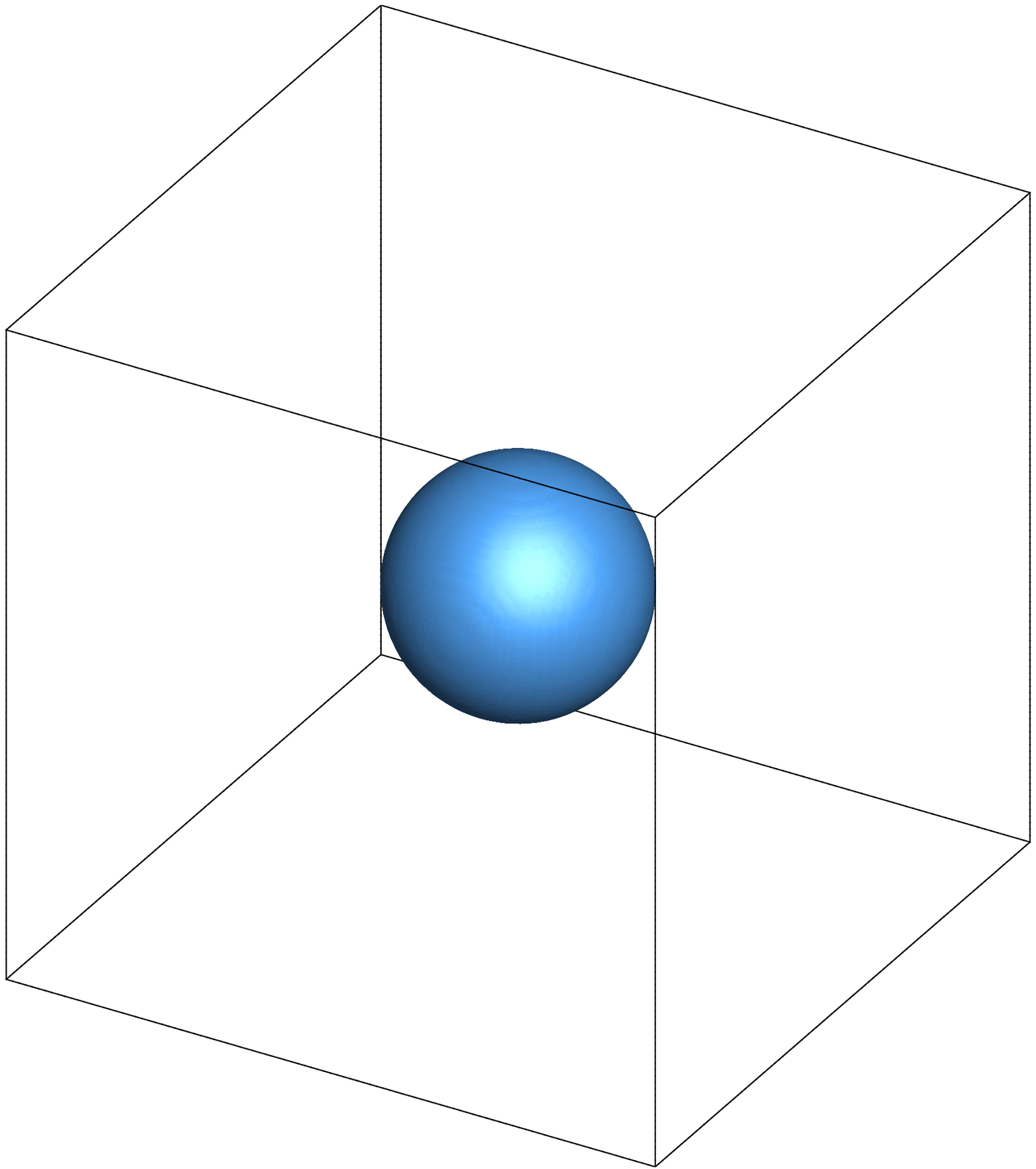}}
	\caption*{ (b) t = 80000$\delta t$}		
	
	\caption{ Snapshots of the density contours obtained by the FDM (left) and the present LB model (right) at (a) t = 40000$\Delta t$ and (b) t = 80000$\Delta t$.      }
	\label{fig5}
\end{figure}

\begin{figure}
	\centering
	\includegraphics[width=0.5\textwidth]{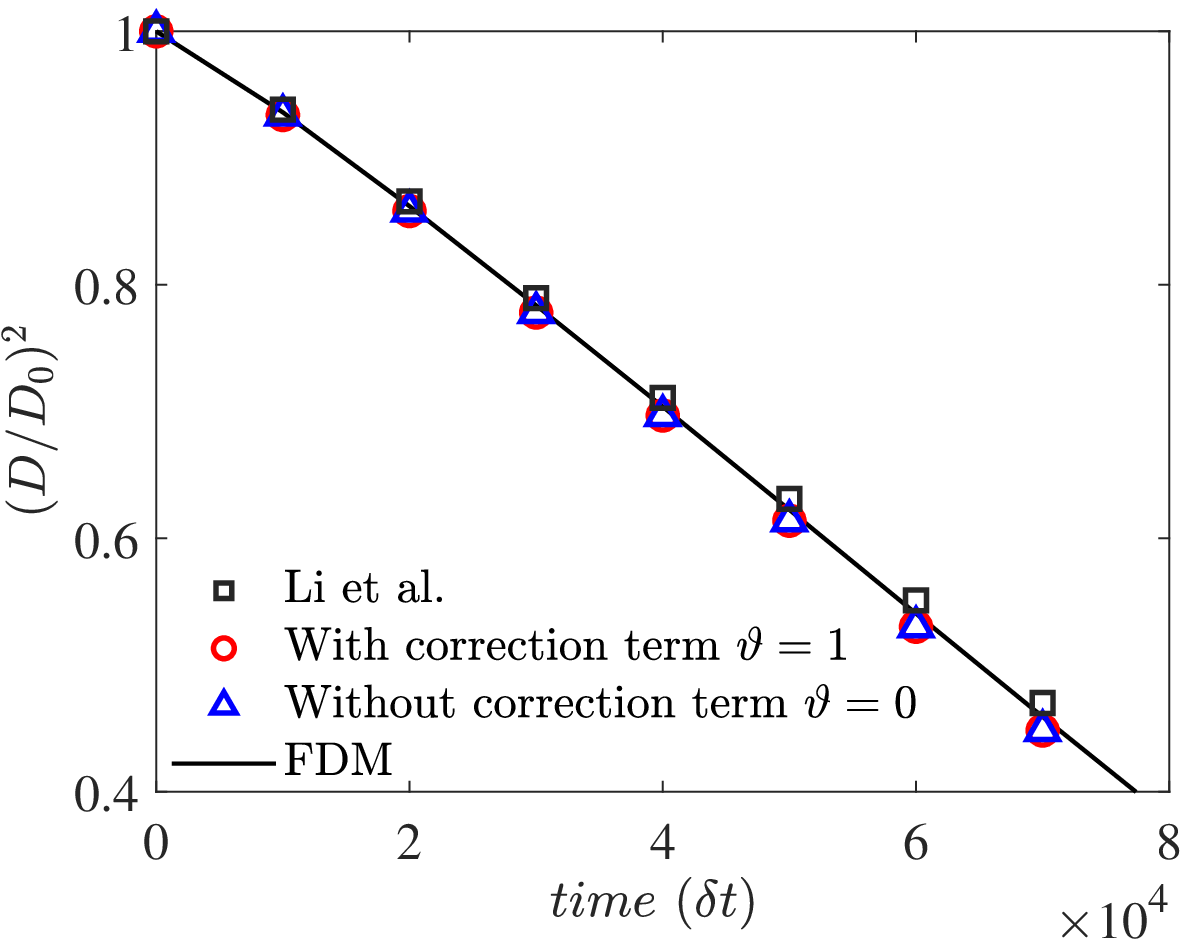}\\		
	\caption{ Time evolution of the square of the dimensionless diameter for different models. }
	\label{fig6}
\end{figure}

Fig. \ref{fig5} presents the distribution of the density contours at two different times, in which the numerical results obtained by using the FDM are also involved for comparison. It can be seen that the droplet tends to evaporate because of the presence of the temperature gradient at the liquid-vapor interface. Also, we note that the difference between the present results and FDM's results is insignificant. To give a quantitative analysis, Fig. \ref{fig6} further gives the variation of the square of droplet diameter with time, where the $d^2$-law is well predicted. Moreover, although the numerical results obtained by the proposed model and Li et al.'s model \cite{LiPRE2022} fit well with the FDM, the present model is more efficient since it is not necessary to calculate the density gradient (see Table \ref{label1}). 

\begin{table}
\caption{Comparisons of the computational performance of different models for  10,000 iteration steps, the custom C code of each case is run on a personal computer with Intel$^\circledR$ Core$^\text{TM}$ i7-9700 CPU @3.0GHz base frequency and 16.0 GB share memory.}
\centering
\begin{tabular}{cccc}
\hline\hline
             & Li et al. \cite{LiPRE2022} & Present model, $\vartheta = 0$ & Present model, $\vartheta = 1$ \\ \hline
CPU time (s) & 57599.4   & 54461.5                        & 55983.6                        \\ \hline\hline
\end{tabular}
\label{label1}
\end{table}

\subsection{Droplet evaporation on a  heated surface}
Droplet evaporation on the heated surface is a fundamental heat transfer phenomenon that is widely encountered in many engineering applications, such as chip manufacturing \cite{DugasLangmuir2005,MarcinichenAE2012}, and spray cooling \cite{ChengRSER2013,GrissomIJHMT2013}, to name but a few. Because of the complex interfacial dynamics and the heat/mass transfer across the liquid-vapor interface, simulating droplet evaporation on a heated surface poses a great challenge for numerical methods \cite{NakoryakovIJHMT2012}. In this subsection, to illustrate the potential of the present model for liquid-vapor phase change problems, we will adopt the proposed thermal LB model to simulate  the droplet evaporation on a heated surface, and compare the numerical results against the experimental data by Dash et al. \cite{DashPRE2014}. Before addressing the initial step of this physical problem, we first recall one dimensionless number using to define the ratio of sensible latent heat absorbed or released during liquid vapor phase change, i.e., Jakob number, 
\begin{equation}
	Ja = \frac{{{c_p}\left( {{T_w} - {T_{sat}}} \right)}}{{{h_{lv}}}},
\end{equation}
where $T_w$ is the wall temperature, and $h_{lv}$ is the specific latent heat.  

According to the experimental \cite{DashPRE2014}, a small droplet with the size of $3 \pm 0.1 \mu l$ is gently deposited on a heated surface at time $t=0$, and initially the air temperature around the droplet is $\left( {21 \pm 0.5} \right){}^ \circ C$, while the wall temperature is set to be $40{}^ \circ C$ with the contact angle being ${120^ \circ }$.  Owing to the droplet is so small that the gravity effect is neglected in the experimental. The dimensionless Jakob number for this problem is given as $Ja=0.036$. Our simulations are conducted in a 3D domain with the grid resolution being  $160 \times 160 \times 100$, and initially a droplet of the diameter equaling $70$ lattice is located at the center of the heated wall, and its temperature is given by $T_s=0.86T_c$. As for the boundary conditions,  the periodic boundary conditions are used in $x$ and $y$ directions, while the constant-pressure boundary condition is used for the top boundary. In addition, the temperature of the bottom and top boundaries are set to be $T_w$ and $T_s$ in simulations. To determine the wall temperature $T_w$, we first calculate the specific latent heat by following the work of Gong and Cheng \cite{GongIJHMT2013}, i.e, $h_{lv}=0.5032$, then the wall temperature is set as $T_w=0.8931T_c$ under the condition of $c_v=5.0$.

\begin{figure}[H]
	\centering
	\includegraphics[width=0.6\textwidth]{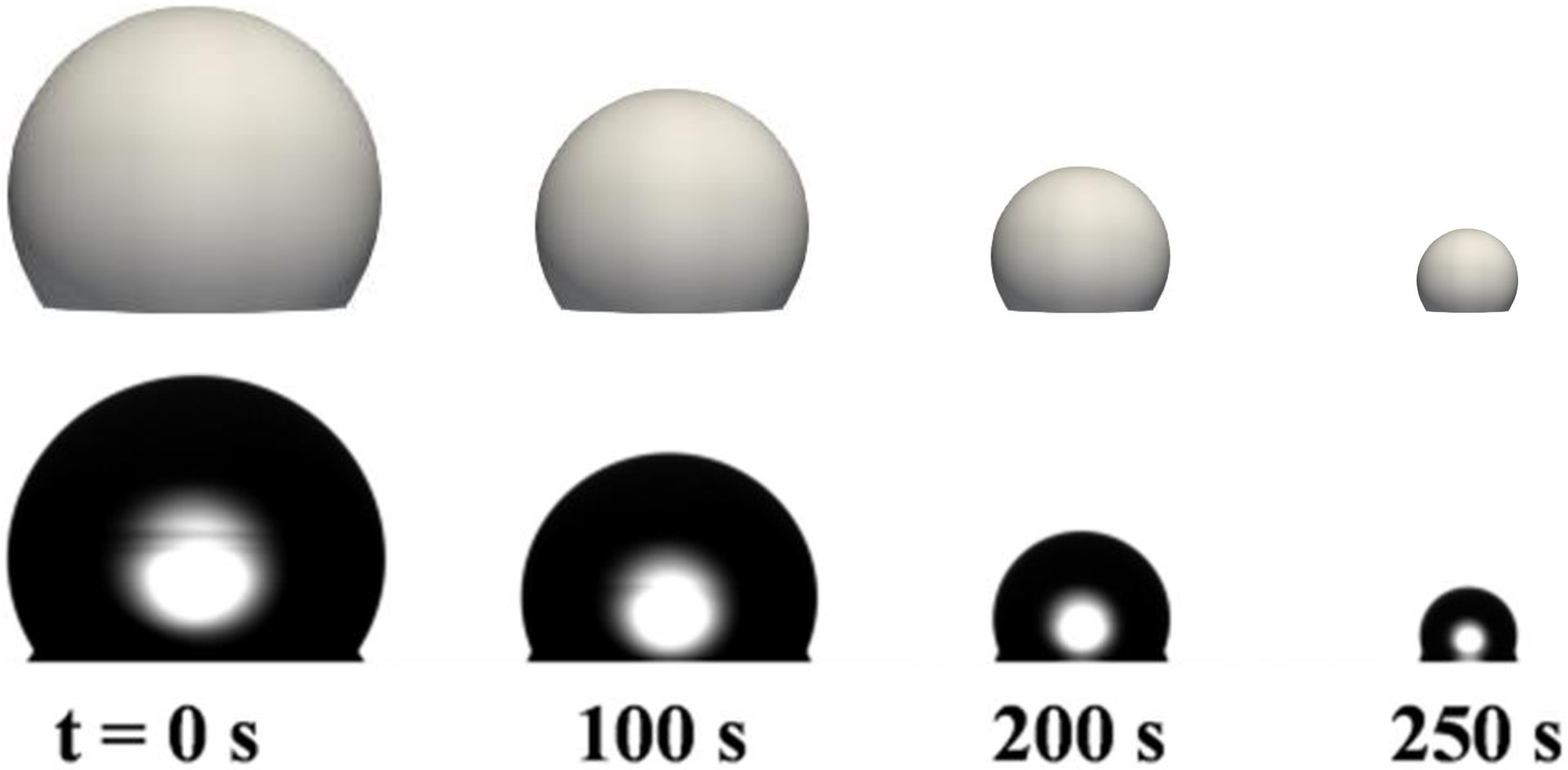}\\		
	\caption{Comparisons of the density contours obtained by the present model (top row) and the experimental (bottom row) \cite{DashPRE2014}. }
	\label{fig7}
\end{figure}

\begin{figure}[H]
	\centering
	\includegraphics[width=0.5\textwidth]{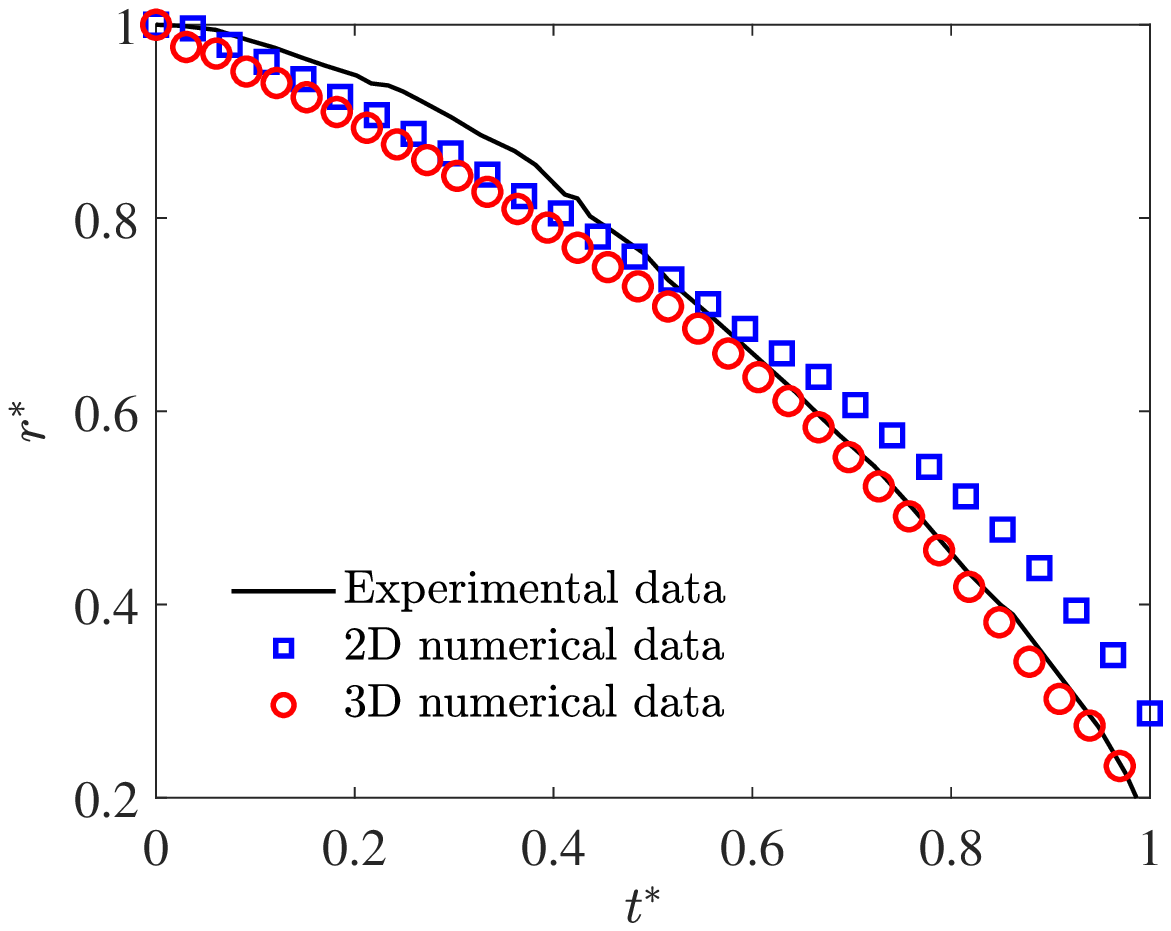}\\		
	\caption{Profiles of the normalized radius $r^*$ with the dimensionless time $t^*$. }
	\label{fig8}
\end{figure}

Fig. \ref{fig7} depicts the interface morphology evolution of the droplet at different times by the proposed model, in which we also provide the experimental results \cite{DashPRE2014}. It is found that due to the different vapor pressure between the liquid-vapor interface and the outlet, the droplet gradually evaporate, and compared to the experimental data, the present model predicts the droplet morphology well.  Fig. \ref{fig8} shows the time evolution of the droplet radius obtained by the current 3D model. For comparisons, the experimental and 2D results are also included. Note that the droplet radius and the time in Fig. \ref{fig8} are normalized when the droplet diameter equals 20 lattices. From Fig. \ref{fig8} we observe a good agreement between our 3D numerical results and experimental data, while evaporation process obtained by 2D model is always slower, which suggests that it is more reasonable to model the liquid-vapor phase change by using a 3D model.        

\section{Conclusions}
As a fundamental phase change problem, liquid-vapor phase change often arises in nature and scientific research. However, the most widely used approach for modeling such flows is still the finite difference method in the LB community. This paper presents an efficient thermal LB model for simulating 3D liquid-vapor phase change. Compared with some existing approaches, the present model is constructed based on the D3Q7 lattice, and it is unnecessary to calculate the volumetric heat capacity gradient, making the current model more efficient and is possible to hold the advantages of the LB method. On the basis of the Zhou-He boundary scheme, a constant-pressure outflow boundary is further proposed to treat the open-boundary condition. The proposed model is first validated by considering the heat transfer in a 3D saturate liquid-vapor system, and found that the obtained numerical results fit well with the theoretical analysis.  Then, it is validated by simulating the droplet evaporation in an open space where the $d^2$-law is well predicted and the numerical results also in consist with the finite-difference-method. At last, to show the model's capability, we simulate the droplet evaporation on a heated surface, and compare the numerical results with the experimental data. It is found that present model is capable in predicting the contact diameter. In summary, our numerical results indicate that the present model is more simple and efficient, and is a promising candidate for simulating 3D liquid-vapor phase change problems.

\section*{CRediT authorship contribution statement}
Jiangxu Huang: Methodology, Validation, Formal analysis, Visualization. Lei Wang: Conceptualization, Methodology, Writing – original draft, Funding acquisition. Kun He:  Methodology, Visualization, Supervision. ChangSheng Huang : Methodology, Supervision.

\section*{Declaration of Competing Interest}
The authors declare that they have no known competing financial
interests or personal relationships that could have appeared to influence
the work reported in this paper

\section*{Ackowledgements}
L.W. is grateful to Dr. Lei Wang from Huazhong University of Science and Technology  , Dr. Anjie Hu from Southwest University of Science and Technology, and Dr. Qin Lou from University of Shanghai for Science and Technology for their help discussions. This work is supported by the Natural Science Foundation of China (Grant No. 12002320).
\appendix
\section{The non-orthogonal transformation matrices  for D3Q19 and D3Q7 lattices model}
The non-orthogonal transformation matrix of the D3Q19 discrete velocity is chosen as 
\begin{equation}
	\centering
	\mathbf{\hat M}=\left[\begin{array}{rrrrrrrrrrrrrrrrrrrrrrrrrrr}
		1 & 1 & 1 & 1 & 1 & 1 & 1 & 1 & 1 & 1 & 1 & 1 & 1 & 1 & 1 & 1 & 1 & 1 & 1 \\
		0 & 1 & -1 & 0 & 0 & 0 & 0 & 1 & -1 & 1 & -1 & 1 & -1 & 1 & -1 & 0 & 0 & 0 & 0 \\
		0 & 0 & 0 & 1 & -1 & 0 & 0 & 1 & -1 & -1 & 1 & 0 & 0 & 0 & 0 & 1 & -1 & 1 & -1 \\
		0 & 0 & 0 & 0 & 0 & 1 & -1 & 0 & 0 & 0 & 0 & 1 & -1 & -1 & 1 & 1 & -1 & -1 & 1 \\
		0 & 1 & 1 & 1 & 1 & 1 & 1 & 2 & 2 & 2 & 2 & 2 & 2 & 2 & 2 & 2 & 2 & 2 & 2 \\
		0 & 2 & 2 & -1 & -1 & -1 & -1 & 1 & 1 & 1 & 1 & 1 & 1 & 1 & 1 & -2 & -2 & -2 & -2 \\
		0 & 0 & 0 & 1 & 1 & -1 & -1 & 1 & 1 & 1 & 1 & -1 & -1 & -1 & -1 & 0 & 0 & 0 & 0 \\
		0 & 0 & 0 & 0 & 0 & 0 & 0 & 1 & 1 & -1 & -1 & 0 & 0 & 0 & 0 & 0 & 0 & 0 & 0 \\
		0 & 0 & 0 & 0 & 0 & 0 & 0 & 0 & 0 & 0 & 0 & 1 & 1 & -1 & -1 & 0 & 0 & 0 & 0 \\
		0 & 0 & 0 & 0 & 0 & 0 & 0 & 0 & 0 & 0 & 0 & 0 & 0 & 0 & 0 & 1 & 1 & -1 & -1 \\
		0 & 0 & 0 & 0 & 0 & 0 & 0 & 1 & -1 & -1 & 1 & 0 & 0 & 0 & 0 & 0 & 0 & 0 & 0 \\
		0 & 0 & 0 & 0 & 0 & 0 & 0 & 1 & -1 & 1 & -1 & 0 & 0 & 0 & 0 & 0 & 0 & 0 & 0 \\
		0 & 0 & 0 & 0 & 0 & 0 & 0 & 0 & 0 & 0 & 0 & 1 & -1 & -1 & 1 & 0 & 0 & 0 & 0 \\
		0 & 0 & 0 & 0 & 0 & 0 & 0 & 0 & 0 & 0 & 0 & 1 & -1 & 1 & -1 & 0 & 0 & 0 & 0 \\
		0 & 0 & 0 & 0 & 0 & 0 & 0 & 0 & 0 & 0 & 0 & 0 & 0 & 0 & 0 & 1 & -1 & -1 & 1 \\
		0 & 0 & 0 & 0 & 0 & 0 & 0 & 0 & 0 & 0 & 0 & 0 & 0 & 0 & 0 & 1 & -1 & 1 & -1 \\
		0 & 0 & 0 & 0 & 0 & 0 & 0 & 1 & 1 & 1 & 1 & 0 & 0 & 0 & 0 & 0 & 0 & 0 & 0 \\
		0 & 0 & 0 & 0 & 0 & 0 & 0 & 0 & 0 & 0 & 0 & 1 & 1 & 1 & 1 & 0 & 0 & 0 & 0 \\
		0 & 0 & 0 & 0 & 0 & 0 & 0 & 0 & 0 & 0 & 0 & 0 & 0 & 0 & 0 & 1 & 1 & 1 & 1
	\end{array}\right] ,
\end{equation}
The  non-orthogonal transformation matrix for D3Q7 lattice is given as 
\begin{equation}
	\centering
   \mathbf{M}=	\left[\begin{array}{ccccccc}
		1 & 1 & 1 & 1 & 1 & 1 & 1 \\
		0 & 1 & -1 & 0 & 0 & 0 & 0 \\
		0 & 0 & 0 & 1 & -1 & 0 & 0 \\
		0 & 0 & 0 & 0 & 0 & 1 & -1 \\
		0 & 1 & 1 & 1 & 1 & 1 & 1 \\
		0 & 1 & 1 & -1 & -1 & 0 & 0 \\
		0 & 1 & 1 & 0 & 0 & -1 & -1
	\end{array}\right].\end{equation}

\end{document}